\documentclass{amsart}

\usepackage{pstricks}
\usepackage{pst-node,pst-tree}
\usepackage [arrow,curve,matrix,tips,frame]{xy}
\usepackage{amsmath}
\usepackage{amsfonts}
\usepackage{latexsym}
\usepackage{euscript}
\usepackage{amscd}

\theoremstyle{plain}
\newtheorem*{thm*}{\sc Theorem}
\newtheorem*{lem*}{\sc Lemma}
\newtheorem*{pro*}{\sc Proposition}
\newtheorem*{cor*}{\sc Corollary}
\theoremstyle{remark}
\newtheorem*{notation*}{\it Notation}
\newtheorem*{exa*}{\it Example}
\newtheorem*{examples*}{\it Examples}
\theoremstyle{remark}
\newtheorem*{definition*}{\sc Definition}
\newtheorem*{definitions*}{\sc Definitions}
\newtheorem*{rem*}{\it Remark}
\newtheorem*{remarks*}{\it Remarks}
\theoremstyle{plain}
\newtheorem{thm}{\sc Theorem}[section]
\newtheorem{lem}[thm]{\sc Lemma}
\newtheorem{pro}[thm]{\sc Proposition}
\newtheorem{cor}[thm]{\sc Corollary}
\theoremstyle{remark}
\newtheorem{exa}[thm]{\it Example}

\newtheorem{example}[thm]{\it Examples}

\newtheorem{definition}[thm]{\sc Definition}
\newtheorem{rem}[thm]{\it Remark}
\newtheorem{definitions}[thm]{\sc Definitions}

\newpsobject{showgrid}{psgrid}{subgriddiv=1,griddots=10,gridlabels=6pt}

\def\AA{{\mathbb A}}

\def\CC{{\mathbb C}}

\def\NN{{\mathbb N}}

\def\PP{{\mathbb P}}
\def\QQ{{\mathbb Q}}
\def\RR{{\mathbb R}}

\def\ZZ{{\mathbb Z}}

\def\Ii{{\mathcal I}}
\def\Jj{{\mathcal J}}

\def\Oo{{\mathcal O}}
\def\Pp{{\mathcal P}}

\def\aaa{{\mathfrak a}}

\def\EEE{{\mathfrak E}}

\def\mmm{{\mathfrak m}}

\def\VVV{{\mathfrak V}}

\def\hbar{{\,\overline{\!h}}}

\def\Kbar{{\,\overline{\!K}}}

\def\wbar{{\,\overline{\!w}}}

\def\Ctilde{{\,\widetilde{\!C}}}

\def\ftilde{\,\widetilde{\!f}}

\def\ptilde{{\,\tilde{p}}}
\def\Ptilde{{\,\widetilde{\!P}}}
\def\qtilde{{\,\tilde{q}}}

\def\Ttilde{{\,\widetilde{\!T}}}

\def\xtilde{{\,\tilde{x}}}
\def\Xtilde{{\,\widetilde{\!X}}}
\def\ytilde{{\,\tilde{y}}}

\def\Ztilde{{\,\widetilde{\!Z}}}

\def\bb{{\boldsymbol{b}}}

\def\ee{{\boldsymbol{e}}}

\def\mm{{\boldsymbol{m}}}

\def\tT{{\boldsymbol{T}}}

\def\vv{{\boldsymbol{v}}}

\def\ww{{\boldsymbol{w}}}

\def\one{\operatorname{\bf 1}}

\def\Bl{\operatorname{Bl}}

\def\except{\operatorname{except}}

\def\lct{\operatorname{lct}}

\def\Newt{\operatorname{Newt}}
\def\Pic{\operatorname{Pic}}

\let\lra=\longrightarrow

\def\ie{{\it i.e.}~}

\def\inv{^{-1}}

\let\phi=\varphi
\let\epsilon=\varepsilon

\newcommand{\floor}[1]{\left\lfloor#1\right\rfloor}

\newcommand{\listspace}{\setlength{\itemsep}{-2pt}}

\title{\sc Log-canonical threshold for curves on a smooth surface}
\author{Marian Aprodu and Daniel Naie}
\thanks{Marian Aprodu was supported in part by CNRS, and by the
ANCS contract CEx05-D11-11/04.10.05. He expresses his
thanks to the LAREMA Angers for hospitality during
the preparation of this work.}
\date{\today}
\address{Marian Aprodu, Institute of Mathematics ``Simion Stoilow'' of
  the Romanian Academy, P.O. BOX 1-764, RO-014700 Bucharest, Romania and
  \c Scoala Normal\u a Superioar\u a, Calea Grivitei 21, RO-010702
  Bucharest, Romania}
\email{{Marian.Aprodu@imar.ro}}
\address{Daniel Naie, Universit\'e d'Angers,
2, Bd Lavoisier, F-49045 Angers, France}
\email{{daniel.naie@univ-angers.fr}}
\begin{document}

\begin{abstract}
It is shown that the log-canonical threshold of a curve with an
isolated singularity is computed by the term ideal of the curve in a
suitable system of local parameters at the singularity.  The proof
uses the Enriques diagram of the singularity and shows that the
log-canonical threshold depends only on a {\it non-degenerate} path of
that diagram.
\end{abstract}

\maketitle

\section{Introduction}

Let $X$ be an affine variety with coordinate ring $A=\CC[X]$.  Any
ideal $\aaa\subset A$ determines a sheaf on $X$ that will also be
denoted by $\aaa\subset\Oo_X$.  Let $f\in A$ with an isolated
singularity at $P$ and let $\aaa_f$ be its term ideal in a system of
local parameters $\Pp$ at $P$.  It is known that for a rational
number $\xi>0$ one can attach to the divisor $(f)$ a collection of
multiplier ideals $\Jj(\xi\cdot(f))$ that starts at $\Oo_X$,
diminishes exactly when $\xi$ equals a {\it jumping number}---they
represent an increasing discrete sequence of rationals---and finally
ends at $\Ii_{(f)}$.  The multiplier ideals reflect the singularity of
the rational divisor $\xi (f)$.  Multiplier ideals and jumping numbers
are also attached to ideal sheaves in a similar way.  For $0<\xi<1$,
one has $\Jj(\xi\cdot (f))\subset\Jj(\xi\cdot\aaa_f)$ with equality if the
coefficients of $f$ are `sufficiently general'.  The definition and
theorem hereafter will make precise what `sufficiently general' means
and state the equality.

The aim of this paper is to show that the log-canonical threshold, \ie
the first jumping number, of a curve $C$ with an isolated singularity
is always computed by the term ideal in a suitable system of local
parameters, regardless the genericity of the curve's equation.

\begin{definition}[see \cite{ArGuVa}]
Let $\Pp$ be a system of local parameters and let $f\in\CC[X]$ with term
ideal $\aaa_f$.  Let $\Newt(\aaa_f)$ be the Newton polygon associated
to $\aaa_f$.  Given any face $\sigma$ of $\Newt(\aaa_f)$ denote by
$f_\sigma$ the sum of those terms of $f$ corresponding to points lying
on $\sigma$. The function $f$ is said to be {\it non-degenerate along
$\sigma$} if and only if the $1$-form $df_\sigma$ is nowhere vanishing
on the torus $(\CC^\ast)^n$.  The function $f$ is said to have {\it
non-degenerate principal part} if this condition holds for every
compact face $\sigma$ of $\Newt(\aaa_f)$.
\end{definition}

\begin{thm}[see {\cite[Proposition 9.2.28]{Ho2}}]
If $f$ has non-degenerate principal part, then for every $0<\xi<1$ 
\[
  \Jj(\xi\cdot(f)) = \Jj(\xi\cdot\aaa_f)
\]
in a neighbourhood of $P$.
\end{thm}

To state the main result we need the following:

\begin{definition}
Let $C$ be a curve with an isolated singularity at $P$ and let $K$ be
the cluster associated to the minimal log resolution of the
singularity of $C$.  A system of local parameters at $P$ is said to be
an {\it adapted system of local parameters at $P$ } for $C$ if the
divisor $(y)$ contains a free point of $K$ and $(x)$ does not.
\end{definition}

\begin{thm*}[\ref{th:theResult}]
If $C$ is a germ of curve with an isolated singularity at $P$, then
\[
  \lct(C;P)
  = \min\,\{\lct(\aaa_{C,\Pp}) \mid \text{$\Pp$ adapted
  system of local parameters at $P$  for $C$}\}
\]
where $\aaa_{C,\Pp}$ is the term ideal of $C$ with respect to $\Pp$.
\end{thm*}

The proof is given in \S\ref{s:last} and is based on the relation
between monomial ideals and Enriques diagrams and on the study of the
log-canonical threshold of unibranch Enriques diagrams.  The former is
established in \S\ref{s:edMonIdeals} and the latter in
\S\ref{s:unibranch} showing that for a unibranch Enriques diagram the
log-canonical threshold is computed at a point belonging to its
non-degenerate part; see Section \ref{s:edMonIdeals} for definitions.

\section{Preliminaries and notation}
\label{s:notConv}
\subsection{Log-canonical threshold}

Let $X$ be an affine variety and let $\aaa\subset\Oo_X$ be a non-zero
ideal sheaf.  A {\it log resolution} of $\aaa$ is a projective
birational mapping $\mu:Y\to X$ with $Y$ non-singular such that
\[
  \mu\inv\aaa := \aaa\cdot\Oo_Y = \Oo_Y(-F)
\]
where $F$ is an effective divisor on $Y$ such that $F+\except(\mu)$
has simple normal crossing support.

\begin{definition*}
Let $\aaa\subset\Oo_X$ be a non-zero ideal sheaf and $\xi>0$ a rational
number.  If $\mu:Y\to X$ is a log resolution of $\aaa$ with
$\aaa\cdot\Oo_Y=\Oo_Y(-F)$, then the multiplier ideal of
$\xi\cdot\aaa$ is 
\[
  \Jj(\xi\cdot\aaa) = \mu_\ast\Oo_Y(K_{Y|X}-\floor{\xi F}).
\]
\end{definition*}

\begin{definition*}
The {\it log-canonical threshold} at $P\in X$ of a non-zero ideal
sheaf is
\[
  \lct(\aaa;P) = \inf\{\xi\in\QQ \mid \Jj(\xi\cdot\aaa)_P\subset\mmm_P\}.
\]
\end{definition*}

The above definitions work for $D\subset X$ a rational effective
divisor.  In case $C$ is a curve and $X$ an affine
smooth surface, if $\mu^\ast C=\sum_\alpha e_\alpha E_\alpha+\Ctilde$ and 
$K_{Y|X}=\sum_\alpha k_\alpha E_\alpha$, then 
\begin{equation} \label{eq:lctCurve}
  \lct(C;P) = \min_\alpha\frac{k_\alpha+1}{e_\alpha},
\end{equation}
\ie $\lct(C;P)$ is the first $\xi>0$ for which each exceptional
divisor appears in $K_{Y|X}-\floor{\mu^\ast(\xi C)}$ with coefficient
$\geq-1$ and at least one has coefficient $=-1$.

\begin{rem*}
Let $\aaa$ or $D$, $\mu:Y\to X$ and $F$ as above.  The multiplier
ideals of $D$ or $\aaa$ do not depend on the log resolution, see
\cite[II. Thm.9.2.11]{La}.  On the contrary, if $\xi$ is non-integral,
the ideal $\mu_\ast\Oo_Y(-\floor{\xi F})$ may depend on the log
resolution.  The simplest example is given by the ideal $(x^2,y^4)$.
The general element of $\aaa$ defines a curve with a tacnode at the
origin in the affine plane.  The minimal log-resolution of the ideal
is given by $\mu$, the composition of two blow-ups.  If $\mu'$ is
$\mu$ composed with the blowing up of the intersection point of the
two exceptional divisors, then 
$\mu_\ast\Oo_Y(-\floor{1/6 F})=\Oo_{\AA^2}$ and
$\mu'_\ast\Oo_{Y'}(-\floor{1/6 F'})=\Ii_O$.  Now, if $\xi=1$, then
$\overline{\aaa}=\mu_\ast\Oo_Y(-F)$ (see~{\cite[II. p. 216-219]{La}}).
\end{rem*}

\begin{rem*}
If $X$ is a smooth affine surface then, for an ideal sheaf or a curve
$C$ with a singularity at $P$, there exists a unique minimal log resolution.
\end{rem*}

\subsection{Term ideals}
Let $x,y$ be regular functions on $X$ that form a system of local
parameters $\Pp$ at $P$.  The functions $x$ and $y$ will be seen as
coordinates around $P$ through the monomorphism
$\Oo_{X,P}\hookrightarrow\CC[[x,y]]$ which associates to each
$f\in\Oo_{X,P}$ its Taylor power series---the injection of $\Oo_{X,P}$
in its formal completion at $\mmm_P$.
With $\Pp$ a system of local parameters at $P$ given, if $f\in\Oo_X$,
the {\it term ideal} of $f$ is the ideal $\aaa_{f,\Pp}$ generated by
the monomials of $f$ with respect to $\Pp$.  If $C=(f)$ the notation
$\aaa_{C,\Pp}$ will also be used for $\aaa_{f,\Pp}$.

\subsection{Monomial ideals, Newton polygons and Howald's result}

The multiplier ideals and the jumping numbers are hard to compute in
general.  An exception is the class of monomial ideals.  Howald
describes in \cite{Ho} their multiplier ideals by a combinatorial
formula.  To state the result let $\aaa\subset \CC[X]$ be
a monomial ideal, \ie an ideal generated by monomials of the form
$x^\mm=x_1^{m_1}\cdots x_n^{m_n}$ in a system of local parameters
$(x_1,\ldots,x_n)$, with $\mm\in\ZZ^n$.  Such a monomial
ideal $\aaa$ can be identified with the set of exponents in $\ZZ^n$ of
its monomials.  The convex hull of this set in $\RR^n=\ZZ^n\otimes\RR$
is called the {\it Newton polyhedron} of $\aaa$ and it is denoted by
$\Newt(\aaa)$.  If $\one=(1,\ldots,1)\in\ZZ^n$, then:

\begin{thm}[Howald]
  \label{th:mIdeals}
Let $\aaa\subset \CC[X]$ be a monomial ideal in a system of local
parameters.  Then for every $\xi>0$, the multiplier ideal
$\Jj(\xi\cdot\aaa)$ is the monomial ideal generated by all monomials
$x^\mm$ with the vector $\mm$ satisfying
\[
   \mm + \one \in (\xi\cdot\Newt(\aaa))^\circ.
\]
\end{thm}

For example if the polyhedron $\Newt(\aaa)$ is cut out in the first
orthant by the inequalities $g_j(\vv)\geq1$ with non-negative rational
coefficients, then 
\begin{equation}
  \label{eq:lctForMIdeals}
  \lct(\aaa) = \min_j g_j(\one).
\end{equation}
The picture of the Newton polygon of the monomial ideal 
$\aaa = (x^8,x^3y^2,y^4)$ in the figure below shows, using Howald's
result, that $\Jj(5/12\cdot\aaa)=(x,y)$.  Note that even though
$(0,0)+\one$ lies in $\Newt(\aaa)$ it does not lie in its interior.
Therefore, the monomial $1$ corresponding to $(0,0)$ does {\it not}
belong to $\Jj(5/12\cdot\aaa)$.  But for all $c<5/12$,
$1\in\Jj(\xi\cdot\aaa)$.  Hence $\lct(\aaa)=5/12$.

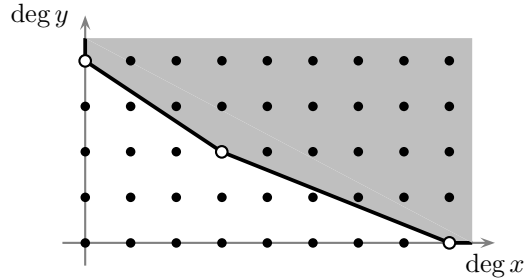
\begin{center}
\begin{pspicture}(0,0)(7,4)%\showgrid
  \psset{unit=4ex}
  \rput(1,1){
    \psline[linecolor=gray,arrowsize=4pt]{->}(-.5,0)(9,0)
    \rput(9,-.5){$\deg x$}
    \psline[linecolor=gray,arrowsize=4pt]{->}(0,-.5)(0,5)
    \rput(-1,5){$\deg y$}
    \pscustom[fillstyle=solid,fillcolor=lightgray,linewidth=.5mm]{%
      \psline(0,4.5)(0,4)(3,2)(8,0)(8.5,0)}
    \pscustom[fillstyle=solid,fillcolor=lightgray,linestyle=none]{%
      \psline(8.5,0)(8.5,4.5)(0,4.5)}
    \multirput(0,0)(0,1){5}{\multirput(0,0)(1,0){9}{\psdots*(0,0)}}
    \pscircle[fillstyle=solid,fillcolor=white](8,0){.15}
    \pscircle[fillstyle=solid,fillcolor=white](3,2){.15}
    \pscircle[fillstyle=solid,fillcolor=white](0,4){.15}
  }
\end{pspicture}
\end{center}

\begin{exa}
  \label{e:useful}
Let $f(x,y)=(x^3-y^2)^2-x^5y$.  Its term ideal is generated by $x^6$,
$x^5y$, $x^3y^2$ and $y^4$ and its Newton polygon has a unique compact
face $\sigma$ given by $m/6+n/4=1$.  It gives $\lct(\aaa_f)=5/12$.
Notice that The polynomial $f$ does not have a non-degenerate
principal part.
\end{exa}

\subsection{Curves with isolated singularities and clusters}

Let $C$ be the germ of a curve on a smooth surface having an isolated
singularity at $P$.  If $\mu:Y\to X$ is the minimal log resolution of
$C$ at $P$,
\[
  \mu^\ast C = \Ctilde + D = \Ctilde + \sum_\alpha e_\alpha E_\alpha.
\]
The proof of Theorem \ref{th:theResult} will mainly deal with the
configuration formed by the strict transforms of the exceptional
divisors $E_\alpha$.  To prepare the way for the proof we need to
formalize the setup and recall some results from the theory of
clusters.

Let $Y=Y_{r+1}\to Y_r\to\cdots\to Y_1=X$ be the decomposition of
$\mu:Y\to X$ into successive blowing ups with
$Y_{\alpha+1}=\Bl_{P_\alpha}Y_\alpha$.  Each point $P_\alpha$ is
infinitely near to $P=P_1$ and has an associated exceptional divisor
on $Y_{\alpha+1}$.  The exceptional divisor and all its strict
transforms will be denoted by $E_\alpha$.  Its total transforms will
be denoted by $W_\alpha$.  When needed, $E_\alpha^{(\beta)}$ will
design the strict transform of $E_\alpha$ on $Y_\beta$ and similarly
for the total transform.  For example
$W_\alpha^{(\alpha+1)}=E_\alpha^{(\alpha+1)}$.  The strict transforms
$E_\alpha$ and the the total transforms $W_\alpha$ form two bases of
the $\ZZ$-module 
$\Lambda_C=\bigoplus_\alpha\ZZ E_\alpha\subset\Pic Y$.  The divisor
$D$ becomes $D=\sum_\alpha w_\alpha W_\alpha$ in the basis of the
$W$'s.

\begin{definition*}
A {\it cluster} in $X$ centered at a smooth point $P$ is a finite set
of weighted infinitely near points to $P$,
$K=\{P_1^{w_1},\ldots,P_r^{w_r}\}$, with $P_1=P$ and such that the
ordering of the points is compatible with the partial order of the
infinitely near points%
\footnote{If $Q$ and $R$ are infinitely near to $P$, then the point
  $Q$ precedes the point $R$ if and only if $R$ is infinitely near to
  $Q$.}.  
The point $P_1$ is called the {\it proper} point of the cluster.
\end{definition*}

The combinatorics of the configuration of the strict transforms on $Y$
is encoded in the notion of {\it proximity} for the points of the
cluster: a point $P_\beta$ is said to be {\it proximate} to
$P_\alpha$, $P_\beta\prec P_\alpha$, if $P_\beta$ lies on
$E_\alpha^{(\beta)}\subset Y_\beta$.  Besides, a point that is
infinitely near, \ie that is not proper, is always proximate to at
most two other points of the cluster.  It is said to be {\it free} if
it is proximate to exactly one point and {\it satellite} if it is
proximate to exactly two points of the cluster.

Let $\Pi=||p_{\alpha\beta}||$ be the decomposition matrix of the
strict transforms in terms of the total transforms on $Y$; it is also
called the proximity matrix of the cluster.  Since
\begin{equation} \label{eq:defOfPi}
  E_\alpha=W_\alpha-\sum_{P_\beta\prec P_\alpha}W_\beta,
\end{equation}
$p_{\alpha\alpha}=1$ for any $\alpha$ and $p_{\alpha\beta}$ equals
$-1$ if $P_\beta$ is proximate to $P_\alpha$ and $0$ if not.  Along
the $\alpha$ column of $\Pi$ the non-zero elements not on the diagonal
correspond to the points to which $P_\alpha$ is a satellite.  The
matrix $-\Pi\cdot^t\!\Pi$ is the intersection matrix of the curves
$E_\alpha$ on the surface $Y$.  Since the intersection matrix of the
curves $W_\alpha$ is minus the identity, there exists effective
divisors $B_\alpha$ that form the dual basis for the divisors
$-E_\alpha$ with respect to the intersection form.  In the sequel this
basis will be referred to as the {\it branch basis}.  The
decomposition matrix of the basis of strict transforms in terms of the
branch basis is $\Pi\cdot^t\!\Pi$.

\subsection{Unloaded clusters}
 
Let $K=\{P_1^{w_1},\ldots,P_r^{w_r}\}$ be a cluster centered at $P$.
As we have noticed, it defines a divisor $D_K=\sum w_\alpha W_\alpha$
on $Y$ and an ideal sheaf $\mu_\ast\Oo_Y(-D_K)$ on $X$, or
equivalently a subscheme $Z_K$ of $X$.  A $0$-dimensional subscheme
defined by a cluster will be called a {\it complete subscheme}.  The
lemma below clarifies the relation between the divisor $D_K$ and the
ideal sheaf $\mu_\ast\Oo_Y(-D_K)$ or equivalently, the complete
subscheme $Z_K$.

\begin{lem}
Let $D_K=\sum_\alpha b_\alpha B_\alpha$.  If $b_\beta<0$ for a certain
$\beta$, then
\[
  \mu_\ast\Oo_Y(-D_K) = \mu_\ast\Oo_Y(-D_K-E_\beta).
\]
\end{lem}

\proof
We take $\mu_\ast$ on the exact sequence 
\[
  0 \lra 
  \Oo_Y(-D_K-E_\beta) \lra
  \Oo_Y(-D_K) \lra
  \Oo_{E_\beta}(-D_K\mid_{E_\beta}) \lra 0.
\]
Since 
$\deg(-D_K\mid_{E_\beta})=
  -(\sum b_\alpha B_\alpha)\cdot E_\beta=b_\beta<0$, we have that 
$$
\mu_\ast\Oo_{E_\beta}(-D_K\mid_{E_\beta})=0.
$$
\qed

\begin{definition*}
A cluster $K$ is said to satisfy the {\it proximity relations} if for
every $P_\alpha$ in $K$, 
\[
  \wbar_\alpha=\sum_{P_\beta\prec P_\alpha}w_\beta \leq w_\alpha.
\]  
\end{definition*}

\begin{cor}[see \cite{Ca}, Theorem 4.2]
  \label{l:unloadingProcedure} 
Let $K=\{P_1^{w_1},\ldots,P_r^{w_r}\}$ be a cluster that
contains a point $P_\alpha$ at which the proximity relation is not
satisfied.  If $K'=\{P_1^{w'_1},\ldots,P_r^{w'_r}\}$ is the
cluster defined by $w'_\alpha=w_\alpha+1$, $w'_\beta=w_\beta-1$ for
every $\beta$ with $P_\beta$ proximate to $P_\alpha$, and
$w'_\gamma=w_\gamma$ otherwise, then $K$ and $K'$ define the same
subscheme in $X$, \ie $\mu_\ast\Oo_Y(-D_K)=\mu_\ast\Oo_Y(-D_{K'})$. 
\end{cor}

\proof
Let 
$D_K=\sum_\alpha w_\alpha W_\alpha=\sum e_\alpha
E_\alpha=\sum_\alpha b_\alpha B_\alpha$ 
and $D_{K'}=\sum_\alpha b'_\alpha B_\alpha$.  The coefficients
$b_\alpha$ are given by
$\bb=\ee\cdot\Pi\cdot^t\!\Pi=\ww\cdot^t\!\Pi=\ww-\overline{\ww}$,
with $\bb=(b_1,\ldots,b_r)$ and so on.  Then
\[
  \bb' = \ww'-\overline{\ww}' 
  = \ww-\overline{\ww}+(\Pi\cdot^t\!\Pi)_\alpha
  =\bb+(\Pi\cdot^t\!\Pi)_\alpha
\]
and
\[
  \ee' = \bb'\cdot(\Pi\cdot^t\!\Pi)^{-1}
  = \bb\cdot(\Pi\cdot^t\!\Pi)^{-1}+
  (\Pi\cdot^t\!\Pi)_\alpha\cdot(\Pi\cdot^t\!\Pi)^{-1} 
  = \ee +(0,\ldots,1,\ldots,0),
\]
hence $D_{K'}=D_K+E_\alpha$.  But $b_\alpha=w_\alpha-\wbar_\alpha<0$
and the result follows form the previous lemma.  
\qed

The cluster $K'$ in the above corollary is said to be obtained from
$K$ by the unloading procedure.  Starting from $K$, iterated
applications of this procedure lead to a cluster $\Kbar$ that
satisfies the proximity relations and defines the same subscheme in
$X$.  The cluster $\Kbar$ is called the unloaded cluster associated to
$K$.  So a cluster and the unloaded cluster associated to it defines
the same complete subscheme of $X$.

\begin{cor}
  \label{c:unloaded}
A cluster is unloaded if and only if the coefficients of its divisor
in the branch basis are non-negative.
\end{cor}

\begin{exa*}
Let $\{P_1^5,P_2^2,P_3^2,P_4^1,P_5^1\}$ and 
$\{P_1^4,P_2^2,P_3^0,P_4^2,P_5^1\}$ be two clusters with the proximity
matrix
\[
\begin{bmatrix}
  1 & -1 & -1 & -1 &   \\
    & 1  & -1 &    &   \\
    &    & 1  & -1 & -1\\
    &    &    & 1  & -1\\
    &&&&1
\end{bmatrix}.
\]
The former is unloaded; it encodes the log-resolution of the
singularity $x^5-y^7=0$, \ie its associated ideal is the integral
closure of $(x^5,y^7)$.  The latter does not satisfy the proximity
relation at $P_3$.  The unloaded associated cluster is
$\{P_1^4,P_2^2,P_3^1,P_4^1,P_5^0\}$ and the associated divisor
$B_2+B_4$.
\end{exa*}

\subsection{Clusters and Enriques diagrams}

For our purposes, it will be convenient to work with Enriques diagrams
instead of clusters.  The points of a cluster $K$, their weights and
proximity relations were encoded by Enriques in an appropriate tree
diagram now called the {\it Enriques diagram} of the cluster (see
\cite{Ca, EnCh, Ev}).  If the weights are omitted, the tree reflects
the combinatorics of the configuration of the strict transforms
$E_\alpha\subset Y$.

\begin{definition*}[see \cite{EnCh, Ca, Ev}]
An {\it Enriques tree} is a couple $(T,\epsilon_T)$, where $T=T(\VVV,\EEE)$
is an oriented tree (a graph without loops) with a single {\it root},
with $\VVV$ the set of vertices and $\EEE$ the set of edges, and where
$\epsilon_T$ is a map
\[
  \epsilon_T:\EEE\to\{\text{`slant'}, \text{`horizontal'},\text{`vertical'}\}.
\] 
fixing the graphical representation of the edges.
An {\it Enriques diagram} is an weighted or labeled Enriques tree. 
\end{definition*}

\begin{definition*}
Let $T$ be an Enriques tree.  A horizontal (respectively vertical)
{\it $L$-shaped branch} of $T$ is a path of length $\geq1$ such that
all edges, but the first, are horizontal (respectively vertical)
through $\epsilon_T$.  An $L$-shaped branch is proper if it contains at
least two edges.  A {\it maximal $L$-shaped branch} is an $L$-shaped
branch that cannot be continued to a longer one.  
\end{definition*}

\begin{lem}[see {\cite[Proposition 1.2]{Ev}}]
  \label{l:equivalence}
There exists a unique map from the set of clusters in $X$ centered at
a smooth point $P$ to the set of Enriques diagrams such that:
\begin{itemize}\setlength{\itemsep}{-.25\baselineskip}
\item 
For every cluster $K=\{P_1^{w_1},\ldots,P_r^{w_r}\}$ the set of
vertices of the image tree is $\VVV=\{P_1,\ldots,P_r\}$ with the
weights given by the integers $w_1,w_2,\ldots,w_r$.  The root of the
tree is the proper point.
\item
At every point ends at most one edge.
\item
A point $P_\alpha$ is satellite if and only if there is either a
horizontal or a vertical edge that ends at the vertex $P_\alpha$.
\item
If there is an edge that begins at the vertex $P_\alpha$ and ends at
the vertex $P_\beta$ then $P_\beta\in E_\alpha^{(\beta)}$, and the
converse is true if $P_\beta$ is free.
\item
The point $P_\beta$ is proximate to $P_\alpha$ if and only if there is
an $L$-shaped branch that starts at $P_\alpha$ and ends at $P_\beta$.
\item
The strict transforms $E_\alpha$ and $E_\beta$ intersect on $Y$ if and
only if the Enriques diagram contains a maximal $L$-shaped branch that
has $P_\alpha$ and $P_\beta$ as its extremities. 
\item 
An edge that begins at a vertex of a free point and ends at a vertex of
a satellite point is horizontal. 
\end{itemize}
\end{lem}

Due to this correspondence, in the sequel, we shall freely talk about
the vertices of an Enriques tree as being free or satellite.  Now, if
$\tT$ is an Enriques diagram and $K$ a corresponding cluster, the
divisor $D_K$ will also be denoted by $D_\tT$ and will be considered
as {\it associated divisor} to $\tT$. Moreover, if the cluster is
unloaded, the log-canonical threshold of the cluster is defined using
(\ref{eq:lctCurve}) and the notation $\lct(\tT)$ and $\lct(D_K)$ will
be used.

\begin{example} \label{ex:tpq}
Let $p<q$ be relatively prime positive integers.  The Enriques tree
$T_{p,q}$ and diagram $\tT_{p,q}$ encode the minimal log resolution of
the curve $x^p-y^q=0$.  More precisely, if $r_0=a_1r_1+r_2$, \ldots,
$r_{m-2}=a_{m-1}r_{m-1}+r_{m}$ and $r_{m-1}=a_mr_m$, with $r_0=q$ and
$r_1=p$, the oriented tree has 
$\VVV=\{P_\alpha \mid 1\leq\alpha\leq a_1+\cdots+a_m\}$ and
$\EEE=\{[P_\alpha P_{\alpha+1}]\mid 1\leq\alpha\leq
a_1+\cdots+a_m-1\}$.  
The map $\epsilon$ is locally constant on the $a_j$ edges 
$[P_\alpha P_{\alpha+1}]$ with 
$a_1+\cdots+a_{j-1}+1\leq\alpha\leq a_1+\cdots+a_j$.  The first
constant value of $\epsilon$---on the first $a_1$ edges---is `slant'.
The other constant values are alternatively either `horizontal' or
`vertical', starting with `horizontal'.  The Enriques trees
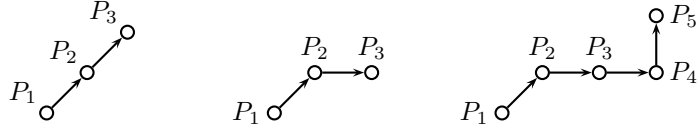
\begin{center}
\begin{pspicture}(0,0)(10,2.25)%\showgrid
  \psset{arrows=->,radius=.1,unit=5ex}
\rput(1,1){
    \Cnode(-.71,-.71){1}
    \Cnode(0,0){2}
    \Cnode(.71,.71){3}
    \ncline{1}{2}
    \ncline{2}{3}
    \uput{\labelsep}[ul](-.71,-.71){$P_1$}
    \uput{\labelsep}[ul](0,0){$P_2$}
    \uput{\labelsep}[ul](.71,.71){$P_3$}}
  \rput(5,1){
    \Cnode(-.71,-.71){1}
    \Cnode(0,0){2}
    \Cnode(1,0){3}
    \ncline{1}{2}
    \ncline{2}{3}
    \uput{\labelsep}[l](-.71,-.71){$P_1$}
    \uput{\labelsep}[u](0,0){$P_2$}
    \uput{\labelsep}[u](1,0){$P_3$}}
  \rput(9,1){
    \Cnode(-.71,-.71){1}
    \Cnode(0,0){2}
    \Cnode(1,0){3}
    \Cnode(2,0){4}
    \Cnode(2,1){5}
    \ncline{1}{2}
    \ncline{2}{3}
    \ncline{3}{4}
    \ncline{4}{5}
    \uput{\labelsep}[l](-.71,-.71){$P_1$}
    \uput{\labelsep}[u](0,0){$P_2$}
    \uput{\labelsep}[u](1,0){$P_3$}
    \uput{\labelsep}[r](2,0){$P_4$}
    \uput{\labelsep}[r](2,1){$P_5$}}
\end{pspicture} 
\end{center}
represent $T_{1,3}$, $T_{2,3}$ and respectively $T_{5,7}$.
The tree $T_{5,7}$ together with the weights $5,2,2,1,1$ and
become the Enriques diagram $\tT_{5,7}$ of the log resolution for
$x^5-y^7=0$.  In general, the weights of $\tT_{p,q}$ are the
corresponding remainders of the Euclidean algorithm.
\end{example}

\begin{exa*}[\ref{e:useful} once more]
The minimal log resolution of $(x^3-y^2)^2-x^5y=0$ needs five blowing
ups with the following Enriques diagram.
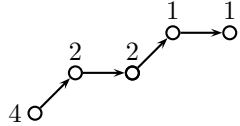
\begin{center}
\begin{pspicture}(0,0)(3,2.25)%\showgrid
  \psset{arrows=->,radius=.1,unit=5ex}
  \rput(1,1){
    \Cnode(-.71,-.71){1}
    \Cnode(0,0){2}
    \Cnode(1,0){3}
    \ncline{1}{2}
    \ncline{2}{3}
    \uput{\labelsep}[l](-.71,-.71){$4$}
    \uput{\labelsep}[u](0,0){$2$}
    \uput{\labelsep}[u](1,0){$2$}}
  \rput(2.71,1.71){
    \Cnode(-.71,-.71){1}
    \Cnode(0,0){2}
    \Cnode(1,0){3}
    \ncline{1}{2}
    \ncline{2}{3}
    \uput{\labelsep}[u](0,0){$1$}
    \uput{\labelsep}[u](1,0){$1$}}
\end{pspicture} 
\end{center}
Using the notation that will be introduced in the \ref{s:last}th section,
the Enriques tree is the connected sum $T_{2,3}\#T_{2,3}$.
\end{exa*}

\section{Newton polygons and Enriques diagrams}
\label{s:edMonIdeals}
We have seen in \S \ref{s:notConv} that to any germ of plane
singularity one can associate an Enriques diagram.  The goal of this
Section is to answer the following question: What are the Enriques
diagrams who correspond to monomial ideals of the polynomial ring in
two variables? Monomial ideals are in some sense the `simplest'
ideals.  The main result is probably well-known to the experts, but due to
the lack of suitable reference and for the coherence of the
exposition we shall state and prove it in this section.

\begin{definitions}
An Enriques tree $(T,\epsilon_T)$ is said to be {\it non-degenerate}
if for every free point $P_\alpha$, the path going from the root to
the vertex $P_\alpha$ contains only free points.  An Enriques tree is
called a {\it binary Enriques tree} if it is non-degenerate, the
outdegree of each vertex is $\leq 2$ and any vertex different from the
root cannot have two proximate free vertices.  An Enriques diagram is
said to be non-degenerate or binary if its tree is.
\end{definitions}

For Enriques trees with the degree of the root equal to $1$ the notion
of {\it union} is defined as follows.  Let $T$ and $T'$ be rooted at
$P$ and $P'$ and let $(S,S')$ be the maximal couple of subtrees of $T$
and $T'$ for which there exists an isomorphism of Enriques trees
$\phi:S\to S'$.  The Enriques tree $T\cup T'$ is the disjoint union of
$T$ and $T'$ modulo the isomorphism $\phi$.  The definition extends to
Enriques diagrams; a vertex belonging neither to $S$ nor to $S'$ keeps
its weight, whereas a vertex corresponding to two vertices identified
through $\phi$ inherits the sum of the respective weights.  It is easy
to see that $T\cup T'$ is again a binary Enriques tree with the degree
of the root equal to $1$.  Moreover, if $\tT$ and $\tT'$ are unloaded,
then $\tT\cup\tT'$ is unloaded.

\begin{exa}
The Enriques diagram $\tT_{5,7}\cup\tT_{4,7}\cup\tT_{3,4}$ is shown
below.
\begin{center}
\begin{pspicture}(0,-.5)(13,2.3)%\showgrid
  \psset{arrows=->,radius=.1,unit=5ex}
\rput(1,.5){
    \Cnode(-.71,-.71){1}
    \Cnode(0,0){2}
    \Cnode(1,0){3}
    \Cnode(2,0){4}
    \Cnode(2,1){5}
    \Cnode(1,1){6}
    \Cnode(1,2){7}
    \ncline{1}{2}
    \ncline{2}{3}
    \ncline{3}{4}
    \ncline{4}{5}
    \ncline{3}{6}
    \ncline{6}{7}
    \uput{\labelsep}[r](-.71,-.71){$12$}
    \uput{\labelsep}[u](0,0){$6$}
    \uput{\labelsep}[d](1,0){$4$}
    \uput{\labelsep}[r](2,0){$2$}
    \uput{\labelsep}[r](2,1){$1$}
    \uput{\labelsep}[r](1,1){$1$}
    \uput{\labelsep}[r](1,2){$1$}
    \uput{\labelsep}[u](3,0){$=$}}
\rput(5.5,1){
    \Cnode(-.71,-.71){1}
    \Cnode(0,0){2}
    \Cnode(1,0){3}
    \Cnode(2,0){4}
    \Cnode(2,1){5}
    \ncline{1}{2}
    \ncline{2}{3}
    \ncline{3}{4}
    \ncline{4}{5}
    \uput{\labelsep}[l](-.71,-.71){$5$}
    \uput{\labelsep}[u](0,0){$2$}
    \uput{\labelsep}[u](1,0){$2$}
    \uput{\labelsep}[r](2,0){$1$}
    \uput{\labelsep}[r](2,1){$1$}
    \uput{\labelsep}[u](3,-.5){$\cup$}}
  \rput(9.75,.5){
    \Cnode(-.71,-.71){1}
    \Cnode(0,0){2}
    \Cnode(1,0){3}
    \Cnode(1,1){4}
    \Cnode(1,2){5}
    \ncline{1}{2}
    \ncline{2}{3}
    \ncline{3}{4}
    \ncline{4}{5}
    \uput{\labelsep}[l](-.71,-.71){$4$}
    \uput{\labelsep}[u](0,0){$3$}
    \uput{\labelsep}[r](1,0){$1$}
    \uput{\labelsep}[r](1,1){$1$}
    \uput{\labelsep}[r](1,2){$1$}
    \uput{\labelsep}[u](2,0){$\cup$}}
  \rput(13,1){
    \Cnode(-.71,-.71){1}
    \Cnode(0,0){2}
    \Cnode(1,0){3}
    \Cnode(2,0){4}
    \ncline{1}{2}
    \ncline{2}{3}
    \ncline{3}{4}
    \uput{\labelsep}[l](-.71,-.71){$3$}
    \uput{\labelsep}[u](0,0){$1$}
    \uput{\labelsep}[u](1,0){$1$}
    \uput{\labelsep}[u](2,0){$1$}}
\end{pspicture} 
\end{center}
\end{exa}

\begin{thm}
  \label{th:polygonDiagram}
Let $\aaa\subset\Oo_X$ be an ideal whose associated subscheme is
$0$-dimensional and centered at $P$.  There exists a system of local
parameters at $P$ such that $\aaa$ is monomial if and only if the
Enriques diagram associated to $\aaa$ is a binary Enriques diagram.
In this case the followings hold:
\begin{enumerate}
  \listspace
\item
  the boundary lines of $\Newt(\aaa)$ have slopes $-q_j/p_j$ and
  lengths $d_j$, with $j\in J$, if and only if the two subdiagrams
  with the degree of the root $=1$ into which the Enriques diagram
  decomposes are 
\[
  \tT' = \bigcup_{j,-\frac{q_j}{p_j}\leq-1} \tT_{p_j,q_j}^{d_j}
  \qquad\text{and}\qquad
  \tT'' = \bigcup_{j,-\frac{q_j}{p_j}>-1} \tT_{p_j,q_j}^{d_j};
\]
\item
  if $P'$ and $P''$ are the extremities of the longest paths
  of free points of these two subdiagrams, then $P'$ and $P''$ lie on
  the coordinate curves.
\end{enumerate}
\end{thm}

\proof
The proof proceeds by induction on the height of the root---the
length of the path from the root to the furthest child.  Suppose the
theorem true when the corresponding Enriques diagram has height of the
root equal to $n$.  Let $\aaa$ be an ideal such that the associated
Enriques diagram is binary and the height of the root equals $n+1$.
Let $\sigma:\Xtilde=\Bl_P(X)\to X$ be the blowing up at $P$ and let
$E$ be the exceptional curve.  Then
\[
  \sigma^\ast\aaa = \Ii_{\Ztilde_1\cup\Ztilde_2}\otimes\Oo_{\Xtilde}(-cE)
\]
with $\Ztilde_1$ and $\Ztilde_2$ $0$-dimensional subschemes supported
at $\Ptilde_1$ and $\Ptilde_2$, and
$c\geq\deg\Ztilde_1|_E+\deg\Ztilde_2|_E$.  The Enriques diagrams
associated to $\Ztilde_1$ and $\Ztilde_2$ are binary and the heights
of the root are $\leq n$.  By the induction hypothesis $\Ztilde_1$ and
$\Ztilde_2$ are monomial in suitable systems of local parameters
$\Pp_j$ at the points $\Ptilde_j$, $j=1,2$, and satisfy the claims in
Theorem \ref{th:polygonDiagram}.  Let $\Pp$ be a local system of local
parameters at $P$ such that $\sigma$ is given locally by
\[
\begin{cases}
  x = \xtilde_1\ytilde_1 \\
  y = \ytilde_1
\end{cases}
\qquad\text{and}\qquad
\begin{cases}
  x = \xtilde_2 \\
  y = \xtilde_2\ytilde_2.
\end{cases}
\]
Using the notion of staircase below and the lemma hereafter with $Z_1$
and $Z_2$ the complete subschemes associated to $\Ztilde_1$ and
$\Ztilde_2$, we conclude that $\aaa$ is monomial and that
$\mu_\ast\Oo(-D_\tT)=\overline{\aaa}$ (the integral closure
$\overline{\aaa}$ of a monomial ideal $\aaa$ is the monomial ideal
spanned by all monomial whose exponents lie in $\Newt(\aaa)$).  The
other claims of the theorem follow from the identity on the staircases
and the way the Enriques diagram $\tT$ is constructed from the
diagrams $\tT_j$, $j=1,2$:
\[
  \VVV(\tT) = \VVV(\tT_1)\cup\VVV(\tT_2)\cup\{P\},
  \qquad
  \EEE(\tT) = \EEE(\tT_1)\cup\EEE(\tT_2)\cup\{[P\Ptilde_1],[P\Ptilde_2]\}.
\]
The map $\epsilon_\tT$ inherits the values from the corresponding maps
$\epsilon_{\tT_j}$ except for the slant values of a possible path of
free points in each $T_j$.  More precisely, the possible path in $T_1$
of free points coming from a horizontal staircase in the previous
step, becomes a path of horizontal edges in $T$.  The possible path in $T_2$
of free points coming from a vertical staircase in the previous
step, becomes a path of vertical edges in $T$.
\qed

In the next lemma it will be easier to work with staircases
instead of Newton polygons.  A staircase is a subset of
$\Sigma\subset\NN^2$ such that its complement verifies
$\Sigma^c+\NN^2\subset\Sigma^c$.  If a system of local parameters at
$P\in X$ is given, to a staircase $\Sigma$, a monomial ideal is
associated $I^\Sigma$.  Conversely, if $\aaa$ is a monomial ideal in a
certain system of local parameters then $\Sigma(\aaa)$, or $\Sigma(Z)$
if $Z$ is the associated $0$-dimensional subscheme, will denote the
corresponding staircase.

Clearly $\Newt(\aaa)$ is the convex hull of $\Sigma(\aaa)^c$, but
$\Newt(\aaa)$ and $\Sigma(\aaa)$ encode the same information only for
integrally closed monomial ideals.

A finite staircase may be seen as a non increasing finite support
sequence of integers $(n_j)_j$, its horizontal slices.  If the finite
staircases $\Sigma'$ and $\Sigma''$ are defined by the sequences
$(n'_j)_j$ and $(n''_j)_j$, then the sequence $(n'_j+n''_j)_j$ defines
the {\it horizontal sum} $\Sigma'+_h\Sigma''$.  The {\it vertical sum}
is similarly defined.  We shall denote by $\Sigma_c$ the staircase
$\{(m,n)\mid m+n<c\}$.

\begin{lem}
  \label{l:stairsLemma}
Let $X$, $\Xtilde$, $P$, $\Pp$, $P_j$ and $\Pp_j$, $j=1,2$ as in the
above proof.  Let $Z_j$ be a $0$-dimensional complete subscheme of
$\Xtilde$ supported at $\Ptilde_j$ and monomial in the system of local
parameters $\Pp_j$, for each $j=1,2$.  If the integer $c$ satisfies
$c\geq\deg Z_1|_E+\deg Z_2|_E$, then there exists a unique
$0$-dimensional complete subscheme $Z$ of $X$ supported at $P$ and
monomial in the system of local parameters $\Pp$ such that
\[ 
  \Ii_Z = \sigma_\ast(\Ii_{Z_1\cup Z_2}\otimes\Oo_{\Xtilde}(-cE))
\qquad\text{and}\qquad
  \sigma^\ast\Ii_Z = \Ii_{Z_1\cup Z_2}\otimes\Oo_{\Xtilde}(-cE).
\]
Moreover
\[
  \Sigma(Z) = \Sigma_c +_v \Sigma(Z_1) +_h \Sigma(Z_2).
\]
\end{lem}

\proof
It is easy to see that if $c\geq\deg Z_1|_E$, then
$\Ii_{Z'}=\sigma_\ast(\Ii_{Z_1}\otimes\Oo_{\Xtilde}(-cE))$ is a
monomial ideal with the staircase $\Sigma_c+_v \Sigma(Z_1)$---if the
Newton polygon of $Z_1$ is defined by $m,n\geq0$ and $g_i(m,n)\geq0$,
a finite number of linear inequalities with positive integer
coefficients for $m$ and $n$, then the Newton polygon of $Z'$ is
defined by $m,n\geq0$, $g_i(m,m+n-c)\geq0$ and $m+n-c\geq0$---and
similarly for $\sigma_\ast(\Ii_{Z_2}\otimes\Oo_{\Xtilde}(-cE))$.  It
follows that $\sigma_\ast(\Ii_{Z_1\cup Z_2}\otimes\Oo_{\Xtilde}(-cE))$
is the monomial ideal
\[
  \sigma_\ast(\Ii_{Z_1}\otimes\Oo_{\Xtilde}(-cE))
  \cap \sigma_\ast(\Ii_{Z_2}\otimes\Oo_{\Xtilde}(-cE))
\]
with the corresponding staircase 
$(\Sigma_c +_v \Sigma(Z_1)) \cup (\Sigma_c +_h \Sigma(Z_2))$ and 
\begin{equation} \label{eq:stairs}
  (\Sigma_c +_v \Sigma(Z_1)) \cup (\Sigma_c +_h \Sigma(Z_2))
  = \Sigma_c +_v \Sigma(Z_1) +_h \Sigma(Z_2)
\end{equation}
if and only if $c\geq\deg Z_1|_E+\deg Z_2|_E$.  
To finish the proof, let $Z$ be the $0$-dimensional subscheme defined
by $\sigma_\ast(\Ii_{Z_1\cup Z_2}\otimes\Oo_{\Xtilde}(-cE))$.  It
is sufficient to notice that $\sigma^\ast\Ii_Z\otimes\Oo_\Xtilde(cE)$
defines a $0$-dimensional subscheme if and only if the identity
(\ref{eq:stairs}) holds.
\qed

\section{Unibranch Enriques trees and log-canonical threshold}
\label{s:unibranch}

The next step is to study the log-canonical threshold for unibranch
Enriques trees.  The proof of the main result will be reduced to this
case. 

\begin{definition}
Let $T$ and $T'$ be unibranch Enriques trees with
$\VVV(T)=\{P_1,\ldots,P_r\}$ and $\VVV(T')=\{P'_1,\ldots,P'_{r'}\}$.
The connected sum of $T$ and $T'$ is the Enriques tree $T\#T'$ with
the set of vertices $\VVV(T\#T')=\VVV(T)\cup\VVV(T')/\{P_r=P'_1\}$,
the set of edges $\EEE(T\#T')=\EEE(T)\cup\EEE(T')$ and the map
$\epsilon_{T\#T'}$ defined by $\epsilon_T$ and $\epsilon_{T'}$ through
the natural restrictions.
\end{definition}

% In what follows we shall use the notation: if $T$ is an Enriques tree,
% then $\Lambda_T\subset\Pic Y$ denotes the $\ZZ$-module
% $\bigoplus_\alpha E_\alpha^T$; $E_\alpha^T$'s denotes the strict
% transforms and $e_\alpha^T$'s the elements of the dual basis for
% $\Lambda_T^\ast$; similarly for $W_\alpha^T$ and $w_\alpha^T$ and
% $B_\alpha^T$ and $b_\alpha^T$. 
In what follows, if $T$ is an Enriques tree, then
$\Lambda_T\subset\Pic Y$ denotes the $\ZZ$-module 
$\bigoplus_\alpha E_\alpha^T$ with $E_\alpha^T$'s the strict
transforms, and $(e_\alpha^T)$ denotes the basis for $\Lambda_T^\ast$,
dual of $(E_\alpha^T)$.  Similarly, $(w_\alpha^T)$ is the dual basis
of the basis of total transforms $(W_\alpha^T)$ and $(b_\alpha^T)$ the
dual basis of the branch basis $(B_\alpha^T)$.

\begin{pro}
  \label{p:main}
Let $T$ be a unibranch Enriques tree that contains at least one proper
$L$-shaped branch and has $r$ vertices.  Let $T'$ be the unibranch
tree $T_{p',q'}$ that contains $r'$ vertices with $p'\geq1$ and
$q'\geq2$ relatively prime integers.  Then the branch basis of the
Enriques tree $S=T\#T'$ satisfies
\[
  \frac{e_r^{S}(B_\alpha^{})}{k_r^{S}+1}
  > \frac{e_{r+r'-1}^{S}(B_\alpha^{})}{k_{r+r'-1}^{S}+1}
\]
for any $1\leq\alpha\leq r+r'-1$, where
$k_{\alpha}^S=e_\alpha^S(\omega^S)$ are the coefficients of the
`relative canonical divisor', $\omega^S=\sum_\alpha W_\alpha^S$, in
the strict transform basis.
\end{pro}

As a consequence we get the second technical ingredient needed for
Theorem \ref{th:theResult}.

\begin{cor}
  \label{c:lctEDiagram}
Let $\tT$ be an unloaded Enriques diagram such that its associated
tree $T$ is unibranch and degenerate.  If $\tT'$ is the Enriques
diagram obtained by taking away the last node from $\tT$, then
\[
  \lct(\tT') = \lct(\tT).
\]
\end{cor}

\proof
Let $D_\tT=\sum_{\alpha=1}^r b_\alpha B_\alpha^T$.  Since $\tT$ is
unloaded, by Corollary \ref{c:unloaded}, $b_\alpha^T(D_\tT)\geq0$ for
any $\alpha$.  Then
\[
  \frac{1}{\lct(\tT)} = \frac{1}{\lct(D_\tT)} 
  = \max_{\beta=1}^r \frac{e_\beta^T(D_\tT)}{k_\beta^T+1}
  = \max_{\beta=1}^r
    \frac{\sum_\alpha b_\alpha e_\beta^T(B_\alpha^T)}{k_\beta^T+1}
  = \max_{\beta=1}^{r-1}
    \frac{\sum_\alpha b_\alpha e_\beta^T(B_\alpha^T)}{k_\beta^T+1},
\]
the last equality being given by Proposition \ref{p:main} since $\tT$
is degenerate.  Hence 
\begin{equation}
  \label{eq:almost}
  \frac{1}{\lct(\tT)} 
  = \max_{\beta=1}^{r-1} \frac{e_\beta^T(D_\tT)}{k_\beta^T+1}.
\end{equation}
Now, let $\tT'$ be the Enriques diagram
obtained by taking away the last node from $\tT$.  Then
\[
  D_{\tT'} = \sum_{\alpha=1}^{r-1}w_\alpha W_\alpha^{T'}
\]
and
\begin{equation}
  \label{eq:ending}
  \frac{1}{\lct(\tT')} 
  = \max_{\beta=1}^{r-1} \frac{e_\beta^{T'}(D_{\tT'})}{k_\beta^{T'}+1}
\end{equation}

\paragraph{Claim} $e_\beta^{T'}(D_{\tT'})=e_\beta^T(D_\tT)$ for any
$1\leq\beta\leq r-1$.

Indeed, let $\sigma^\ast:\Lambda_{T'}\to\Lambda_T$ be
the monomorphism given by $W_\alpha^{T'}\mapsto W_\alpha^T$ and
$\sigma_\ast$ be the dual epimorphism.  Then 
$D_\tT=\sigma^\ast D_{\tT'}+w_rW_r^T$ and, for any $1\leq\beta\leq r-1$,
\[
  e_\beta^T(D_\tT) 
  = e_\beta^T(\sigma^\ast D_{\tT'})+w_re_\beta^T(W_r^T)
  = \sigma_\ast e_\beta^T(D_{\tT'}).
\]
Using the identity (\ref{eq:defOfPi}) and denoting by $\Pi_T$ the
decomposition matrix for the strict and total transforms associated to
$T$, we have 
\[
  \Pi_T =
  \begin{bmatrix}
    \Pi_{T'} & \ast \\
    0        & 1
  \end{bmatrix}.
\]
Hence 
\begin{multline}
  \sigma_\ast e_\beta^T = \sigma_\ast(\ww^T\cdot\Pi_T\inv)_\beta
  = \sigma_\ast
  \bigg(\ww^T\cdot\begin{bmatrix}
    \Pi_{T'}\inv & \ast \\
    0            & 1
  \end{bmatrix}{\bigg)\!}_\beta \\
  = \bigg(\begin{bmatrix} w^{T'}_1 &\ldots& w_{r-1}^{T'}& 0 \end{bmatrix}
  \cdot
  \begin{bmatrix}
    \Pi_{T'}\inv & \ast \\
    0            & 1
  \end{bmatrix}{\bigg)\!}_\beta
  = e_\beta^{T'}
\end{multline}
for any $1\leq\beta\leq r-1$.

\bigskip

% Since $k_\beta^{T'}=k_\beta^T$ for any $1\leq\beta\leq r-1$, using the
% claim and the formulae (\ref{eq:ending}) and (\ref{eq:almost}) the
% corollary follows.
% \qed
Finally, since $k_\beta^{T'}=k_\beta^T$ for any $1\leq\beta\leq r-1$,
using the claim and the formulae (\ref{eq:ending}) and
(\ref{eq:almost}) the corollary follows.  \qed

For the proof of the proposition we need four preliminary lemmas.  The
first three deal with the combinatorics of $T_{p,q}$, see \cite{Na}.

\begin{lem}
  \label{l:theSequences}
If $(f_j)_{-1\leq j\leq m}$ and $(\delta_j)_{1\leq j\leq m+1}$ are two
finite sequences defined by 
\[
\begin{aligned}
  f_j &= f_{j-2}+a_j\delta_j, \quad\text{for any } 1\leq j\leq m, \\
  \delta_j &= \delta_{j-2}+a_{j-1}\,f_{j-2}, 
  \quad\text{for any } 2\leq j\leq m+1 \quad 
\end{aligned}
\]
and such that $f_{-1}=f_0=0$ and $\delta_0=\delta_1=1$, then the
remainder $r_j$ in the Euclid algorithm is given by
$-f_{j-1}q+\delta_jp$ if $j$ is odd and $\delta_jq-f_{j-1}p$ if $j$ is
even.
\end{lem}

\proof
Left to the reader.
\qed

\begin{rem}
  \label{r:lastTerms}
If $m$ is odd, then $f_m=q$ and $\delta_{m+1}=p$, and if $m$ is even,
then $f_m=p$ and $\delta_{m+1}=q$.  Indeed, let us suppose that $m$ is
odd.  Then the equalities follow since for any $1\leq j\leq m$ the
integers $f_j$ and $\delta_{j+1}$ are relatively prime and 
$0=r_{m+1}=\delta_{m+1}q-f_mp$.
\end{rem}

\begin{lem}
  \label{l:last_eCoeffForB}
If $(f_j)_{-1\leq j\leq m}$ and $(\delta_j)_{1\leq j\leq m+1}$ are the
finite sequences defined in Lemma \ref{l:theSequences}, then for any
$1\leq j\leq m$ and any $1\leq k\leq a_j$
\[
  e_r^T(B_{a_1+\cdots+a_{j-1}+k}^T)
  = \begin{cases}
    (f_{j-2}+k\,\delta_j)\,p & \text{if $j$ is odd}\\
    (f_{j-2}+k\,\delta_j)\,q & \text{if $j$ is even.}
    \end{cases}
\]
\end{lem}

\proof
The proof proceeds by induction on $j$ and $k$.  It is
clear for $j=1$ and any $k$.  Suppose that $j$ is even, $k<a_j$ and
that $e_r^T(B_{a_1+\cdots+a_{j-1}+k})=(f_{j-2}+k\,\delta_j)\,q$.  We
recall that $B_\alpha^T$ is given by the Enriques diagram for which the
weight of the point $P_\alpha$ is $w_\alpha=1$, the weights of all the
points that do not precede $P_\alpha$ are $0$, and all the others are
computed by imposing equalities in the proximity relations.  Then
\[
  B_{a_1+\cdots+a_{j-1}+k+1}^T = B_{a_1+\cdots+a_{j-1}+k}^T +
  W_{a_1+\cdots+a_{j-1}+k+1}^T + B_{a_1+\cdots+a_{j-1}}^T
\]
and 
\[
\begin{aligned}
  e_r^T(B_{a_1+\cdots+a_{j-1}+k+1}^T) 
  &=(f_{j-2}+k\,\delta_j)\,q + r_j +
   (f_{j-3}+a_{j-1}\,\delta_{j-1})\,p\\
  &=(f_{j-2}+k\,\delta_j)\,q + \delta_jq-f_{j-1}p + f_{j-1}p\\
  &=(f_{j-2}+(k+1)\delta_j)\,q.
\end{aligned}
\]
The argument is similar in all the other cases, \ie when either 
$k=a_j$ or $j$ odd.
\qed

\begin{lem}
  \label{l:first_wCoeffForB}
For any $1\leq j\leq m$ and any $2\leq k\leq$ either $a_j+1$ if $j<m$
or $a_m$ if not,
\[
  w_1^T(B_{a_1+\cdots+a_{j-1}+k}^T) 
  = \begin{cases}
    \delta_{j-1}+k\,f_{j-1} & \text{if $j$ is odd}\\
    f_{j-2}+k\,\delta_j & \text{if $j$ is even.}
    \end{cases}
\]
\end{lem}

\proof
We know the result for the last branch divisor $B_r^T$, where
$r=a_1+\cdots+a_m$; we have $w_1^T(B_r^T)=p$.  For an arbitrary
$B_\alpha^T$ we reduce the computation to this known situation.

Let us suppose that $j$ is even, the argument being similar in the
other case.  Let $\Ttilde$ be the Enriques tree given by the first
$a_1+\cdots+a_{j-1}+k$ vertices of $T$ and let $\ptilde<\qtilde$ be the
relatively prime positive integers such that $\Ttilde=T_{\ptilde,\qtilde}$.
By what we have just noticed,
\[
  w_1^T(B_{a_1+\cdots+a_{j-1}+k}^T) =
  w_1^\Ttilde(B_{a_1+\cdots+a_{j-1}+k}^\Ttilde) = \ptilde.
\]
The sequences of Lemma \ref{l:theSequences} associated to $\ptilde$
and $\qtilde$ are given by $f_{-1},\ldots,f_{j-1}$ and
$\ftilde_j=f_{j-2}+k\delta_j$, and by $\delta_1,\ldots,\delta_j$ and
$\tilde{\delta}_{j+1}=\delta_{j-1}+kf_{j-1}$.  Here the $f_i$'s and
the $\delta_i$'s are the elements of the sequences associated to $p$
and $q$.  We end the proof using the remark \ref{r:lastTerms}.
\qed

\begin{lem}
  \label{l:decompositionForB}
Let $T$ and $T'$ be two Enriques unibranch trees with $r$ and
respectively $r'$ vertices.  For every $1\leq\beta\leq r'$ the
following formula holds for the connected sum $S=T\# T'$:
\[
  B_{r+\beta-1}^{S} = 
  w_1^{T'}(B_\beta^{T'})\,(B_r^{S}-W_r^{S})
    +\sum_{\alpha'=1}^{r'}w_{\alpha'}^{T'}(B_\beta^{T'})\,
    W_{r+\alpha'-1}^{S}.
\]
\end{lem}

\proof
Using that 
$w_{r+\alpha'}^{S}(B_{r+\beta-1}^{S})
  =w_{\alpha'}^{T'}(B_\beta^{T'})$ for any $1\leq\alpha'\leq r'-1$ and
discarding the exponent $S=T\# T'$, we have
\[
\begin{split}
  B_{r+\beta-1}^{} 
  &= w_r^{}(B_{r+\beta-1}^{})\,B_r^{}
  +\sum_{\alpha'=1}^{r'-1}w_{r+\alpha'}^{}(B_{r+\beta-1}^{})\,
  W_{r+\alpha'}^{}\\
  &= w_1^{T'}(B_\beta^{T'})\,B_r^{}
  +\sum_{\alpha'=1}^{r'-1}w_{\alpha'+1}^{T'}(B_\beta^{T'})\,
  W_{r+\alpha'}^{}\\
  &= w_1^{T'}(B_\beta^{T'})\,(B_r^{}-W_r^{})
  +\sum_{\alpha'=1}^{r'}w_{\alpha'}^{T'}(B_\beta^{T'})\,
  W_{r+\alpha'-1}^{}.  
\end{split} 
\]
\qed

\proof[Proof of Proposition \ref{p:main}]
Let $s=r+r'-1$.  Clearly
$k_r^{S}=e_r(\omega^{S})=e_r^T(\omega^T)$.  As for
$k_s^{S}$, we have
\[
\begin{split}
  k_s^{S}
  &= e_s^{}(\omega^{S}) \\
  &= e_s^{}(\sum_{\alpha=1}^rW_\alpha^{})-e_s^{}(W_r^{})
    +e_s^{}(\sum_{\beta=r}^sW_\beta^{})\\
  &= e_{r'}^{T'}(e_r^T(\omega^T)W_1^{T'})-e_{r'}^{T'}(W_1^{T'})
    +e_{r'}^{T'}(\omega^{T'})\\
  &= (e_r^T(\omega^T)-1)p'+(p'+q'-1).
\end{split}
\]
Now, if $1\leq\alpha\leq r$, then
\[
  e_{r+r'-1}^{}(B_\alpha^{}) 
  = e_{r+r'-1}^{}(e_r^T(B_\alpha)\,W_r^{})
  = e_r^T(B_\alpha)\,e_{r'}(W_1^{T'})
  = e_r^T(B_\alpha)\,p'.
\]
The desired inequality is equivalent to $(q'-p')e_r^T(B_\alpha)>0$
which is satisfied.  
If $\alpha=r+\beta-1$ with $\beta\geq2$, then by Lemma
\ref{l:decompositionForB},
\[
  e_r^{}(B_{r+\beta-1}^{})
  = w_1^{T'}(B_\beta^{T'})\,e_r^{}(B_r^{})
  = e_r^T(B_r^T)\,w_1^{T'}(B_\beta^{T'})
\]
and
\[
\begin{split}
  e_{r+r'-1}^{}(B_{r+\beta-1}^{})
  &= w_1^{T'}(B_\beta^{T'})\,e_{r+r'-1}^{}(B_r^{}-W_r^{})
  +\sum_{\alpha'=1}^{r'}w_{\alpha'}^{T'}(B_\beta^{T'})\,
  e_{r+r'-1}^{}(W_{r+\alpha'-1}^{})\\
  &= w_1^{T'}(B_\beta^{T'})\,(e_r^T(B_r^T)-1)\,e_{r'}^{T'}(W_1^{T'})
  +\sum_{\alpha'=1}^{r'}w_{\alpha'}^{T'}(B_\beta^{T'})\,
  e_{r'}^{T'}(W_{\alpha'}^{T'})\\
  &= (e_r^T(B_r^T)-1)p'\,w_1^{T'}(B_\beta^{T'})+e_{r'}^{T'}(B_\beta^{T'}).
\end{split}   
\]
The inequality we want to prove becomes
\begin{equation}
  \label{eq:theInequality}
  \frac{e_r^T(B_r^T)\, w_1^{T'}(B_\beta^{T'})}{e_r^T(\omega^T)+1}
  > \frac{(e_r^T(B_r^T)-1)p'\,w_1^{T'}(B_\beta^{T'})+
    e_{r'}^{T'}(B_\beta^{T'})}{(e_r^T(\omega^T)+1)p'+(q'-p')}
\end{equation}
and we shall use the lemmas \ref{l:last_eCoeffForB} and
\ref{l:first_wCoeffForB} to establish it.  Let
$\beta=a'_1+\cdots+a'_{j-1}+k$ with $2\leq k\leq a'_j+1$.  Let us
suppose that $j$ is odd---when $j$ is even the argument is similar and
simpler.  The integer $k$ is either $\leq a'_j$ or $=a'_j+1$.  In the
former case,
\[
  w_1^{T'}(B_\beta^{T'}) = \delta'_{j-1}+kf'_{j-1}
  \quad\text{and}\quad
  e_{r'}^{T'}(B_\beta^{T'}) = (f'_{j-2}+k\,\delta'_j)\,p'.
\]
Inequality (\ref{eq:theInequality}) is then equivalent to 
\begin{equation}
  \label{eq:equivInequality1}
  e_r^T(B_r^T)\,(\delta'_{j-1}+kf'_{j-1})\,(q'-p') > 
  (e_r^T(\omega^T)+1)\,(-(\delta'_{j-1}+kf'_{j-1})\,p'
    +(f'_{j-2}+k\,\delta'_j)\,p').
\end{equation}
Since $T$ contains at least one proper $L$-shaped branch
$e_r^T(B_r^T)>e_r^T(\omega^T)+1$.  So, to see that the inequality
(\ref{eq:equivInequality1}) is true, it is sufficient to show
that
\[
  (\delta'_{j-1}+kf'_{j-1})\,(q'-p') 
  \geq  -(\delta'_{j-1}+kf'_{j-1})\,p'+(f'_{j-2}+k\,\delta'_j)\,p',
\]
or equivalently, that 
\[
  \delta'_{j-1}\,q'-f'_{j-2}\,p' \geq k(\delta'_j\,p'-f'_{j-1}q').
\]
By Lemma \ref{l:theSequences}, this last inequality becomes
$r'_{j-1}\geq kr_j$ which is true since $r'_{j-1}=a'_jr'_j+r'_{j+1}$. 
In the latter case, \ie if $k=a'_j+1$, we have 
\[
  w_1^{T'}(B_\beta^{T'}) = \delta'_{j-1}+(a'_j+1)\,f'_{j-1}
  = \delta'_{j+1}+f'_{j-1}
\]
and
\[
  e_{r'}^{T'}(B_\beta^{T'}) = (\delta'_{j+1}+f'_{j-1})\,q'.
\]
The inequality (\ref{eq:theInequality}) becomes equivalent to 
\[
   \frac{e_r^T(B_r^T)}{e_r^T(\omega^T)+1}
  > \frac{e_r^T(B_r^T)p'+(q'-p')}{(e_r^T(\omega^T)+1)p'+(q'-p')}, 
\]
or again to 
\[
  (e_r^T(B_r^T)-e_r^T(\omega^T)-1)(q'-p') > 0.
\]
The first term is positive since $T$ contains at least a proper
$L$-shaped branch.
\qed

\section{The result}
\label{s:last}

The proof of the main result relies on the preceding investigations
where we have established the relationship between the Enriques
diagrams and the monomial ideals and computed the log-canonical
threshold of unibranch Enriques diagrams.

\begin{thm}
  \label{th:theResult}
If $C$ is a germ of curve with an isolated singularity at $P$, then
\[
  \lct(C;P)
  = \min\,\{\lct(\aaa_{C,\Pp}) \mid \text{$\Pp$ adapted
  system of local parameters at $P$  for $C$}\}
\]
where $\aaa_{C,\Pp}$ is the term ideal of $C$ with respect to $\Pp$.
\end{thm}

\proof
We start with a remark.  Let $\tT$ be the Enriques diagram
corresponding to the minimal resolution of $C$ and let $\Pp=(x,y)$ be
an adapted system of parameters for $C$.  Let $\tT(\Pp)$ be the
subdiagram of $\tT$ determined as follows: If $P_\rho$ is the highest
free point on the coordinate curve $(y)$, then $\tT(\Pp)$ is the
biggest binary Enriques subdiagram of $\tT$ that contains as free
points only the ones of the path from $P_0$ to $P_\rho$.  Then
\begin{equation}
  \label{eq:icEd}
  \overline{\aaa}_{C,\Pp} = \mu_\ast\Oo_Y(-D_{\tT(\Pp)}).
\end{equation}
Indeed, $\overline{\aaa}_{C,\Pp}$ is the smallest integrally closed
monomial ideal containing $f$.  It follows that
$\overline{\aaa}_{C,\Pp} \subset \mu_\ast\Oo_Y(-D_{\tT(\Pp)})$ since
$\mu_\ast\Oo_Y(-D_{\tT(\Pp)})$ is monomial by Theorem
\ref{th:polygonDiagram}.  So 
$\tT(\Pp) \subset \tT_{\overline{\aaa}_{C,\Pp}} \subset \tT$.  By
construction, the $\tT(\Pp)$ is the biggest binary Enriques subdiagram
of $\tT$ that contains as free points those between $P_0$ and
$P_\rho$, hence $\tT(\Pp)=\tT_{\overline{\aaa}_{C,\Pp}}$ and the
assertion follows.

Second, let $s(\tT)$ be a subdiagram of $\tT$ whose Enriques tree is a
path from the root to a leaf and that contains a vertex on which $\tT$
realizes its log-canonical threshold.  By Corollary \ref{c:lctEDiagram}
and descending induction 
\[
  \lct(s(\tT))=\lct(s(\tT)^\circ),
\]
hence
\begin{equation}
  \label{eq:stemAndAll}
  \lct(\tT) = \lct(s(\tT))=\lct(s(\tT)^\circ).
\end{equation}
Here, the biggest non-degenerate subdiagram of $\tT$ is denoted by
$\tT^\circ$. 

To finish the proof, we consider for $\Pp$ a system of local
parameters such that one coordinate contains the highest free point of
$s(\tT)^\circ$.  By (\ref{eq:icEd}) and since $s(\tT)^\circ$ is a
subdiagram of $\tT(\Pp)$, 
\[
  \lct(\aaa_{C,\Pp}) = \lct(\overline{\aaa}_{C,\Pp}) 
  = \lct(\mu_\ast\Oo_Y(-D_{\tT(\Pp)}))
  =\lct(\tT(\Pp))
  \leq \lct(s(\tT)^\circ).
\]
Using (\ref{eq:stemAndAll}), 
$\lct(\aaa_{C,\Pp}) \leq \lct(\tT) = \lct(C;P)$
finishing the proof.
\qed

It has to be noticed that if $f$ does not have non-degenerate
principal part, then other jumping numbers may be different for $f$
and for its term ideal as it will be seen below.

\begin{exa*}
[\ref{e:useful} last time]  
For $f(x,y)=(x^3-y^2)^2-x^5y$ we have seen that its Newton polygon is
given by $m/6+n/4-1\geq0$ and that its Enriques diagram is modeled on
$T_{2,3}\#T_{2,3}$ with weights $4,2,2,1,1$.  We have the same
log-canonical threshold for $(f)$ and $\aaa_f$, but the next jumping
number for $\aaa_f$ is $7/12$ and the next jumping number for $f$ is
$15/26<7/12$, showing that the singularity of $(f)$ at the origin is
worse than that defined by a general element of $\aaa_f$ as expected.
\end{exa*}

\end{document}